\documentclass[12pt]{amsart}
\usepackage{amssymb}
\usepackage{isolatin1}    % Diese Zeile kann bei Problemen geloescht werden

\swapnumbers

\newtheorem{thm}{Theorem}[section]
\newtheorem{lem}[thm]{Lemma}

\newtheorem{prop}[thm]{Proposition}

\newtheorem*{prob*}{Open problem}

\theoremstyle{definition}
\newtheorem{defi}[thm]{Definition}

\theoremstyle{remark}
\newtheorem{rem}[thm]{Remark}
\newtheorem*{rem*}{Remark}

\newtheorem{example}[thm]{Example}

\DeclareMathOperator{\s}{span}

\newcommand{\kringel}{\mathbin{\raise1pt\hbox{$\scriptstyle\circ$}}} 
\newcommand{\pkt}{\mathbin{\raise0pt\hbox{$\scriptstyle\bullet$}}}

\newcommand{\C}{\mathbb{C}}

\newcommand{\N}{\mathbb{N}}

\newcommand{\R}{\mathbb{R}}

\newcommand{\ad}{\mathop{\rm ad}}

\newcommand{\Aut}{\mathop{\rm Aut}}
\newcommand{\Der}{\mathop{\rm Der}}
\newcommand{\Pder}{\mathop{\rm Pder}}
\newcommand{\diag}{\mathop{\rm diag}}

\newcommand{\Lg}{\mathfrak{g}}

\newcommand{\Ll}{\mathfrak{l}}

\newcommand{\CA}{\mathcal{A}}

\newcommand{\CM}{\mathcal{M}}

\newcommand{\CF}{\mathcal{F}}
\newcommand{\CI}{\mathcal{I}}

\newcommand{\CT}{\mathcal{T}}

\newcommand{\mi}{\boldsymbol{-}}

\newcommand{\al}{\alpha}
\newcommand{\be}{\beta}
\newcommand{\ga}{\gamma}
\newcommand{\de}{\delta}

\newcommand{\la}{\lambda}
\newcommand{\om}{\omega}
\newcommand{\ze}{\zeta}

\newcommand{\ra}{\rightarrow}  
\renewcommand{\phi}{\varphi}

\begin{document}

\title[Lie algebra prederivations]{Lie algebra prederivations and 
strongly nilpotent Lie algebras} 
%  Die Kurzfassung kommt oben ueber die Seiten, sie steht in eckigen Klammern
%  Auch Autorennamen koennen eine Kurzfassung haben

\author[D. Burde]{Dietrich Burde}
\address{International University of Bremen\\
  School of Science and Engeneering\\
  Campus Ring 1\\
  28759 Bremen\\
  Germany}
\email{d.burde@iu-bremen.de}


\begin{abstract}
We study Lie algebra prederivations. A Lie algebra
admitting a non-singular prederivation is nilpotent. 
We classify filiform Lie algebras admitting a 
non-singular prederivation but no non-singular derivation.
We prove that any $4$-step nilpotent
Lie algebra admits a non-singular prederivation.
\end{abstract}

\maketitle

\section{Introduction}

Let $\Lg$ be a Lie algebra over a field $k$ and $\Der (\Lg)$ its 
derivation algebra. It is a natural question whether $\Lg$ admits
a non-singular derivation or not. Jacobson has proved \cite{JAC} 
that any Lie algebra over a field of characteristic zero admitting 
a non-singular derivation must be nilpotent. He also asked for
the converse, whether any nilpotent Lie algebra admits a non-singular 
derivation. As it turned out, this was not the case. Dixmier and Lister
\cite{DIL} constructed nilpotent Lie algebras possessing only 
nilpotent derivations. They called this class of Lie algebras
{\it characteristically nilpotent} Lie algebras. This class
has been studied extensively later on \cite{LUT}.\\
On the other hand there exist various generalizations of Lie algebra
derivations, see for example \cite{LEL}, \cite{MUL}. 
For so called {\it prederivations} Jacobson's theorem is also true:
any Lie algebra over a field of characteristic zero admitting a 
non-singular prederivation is nilpotent \cite{BAJ}.
Lie algebra prederivations have been studied in connection with 
bi-invariant semi-Riemannian metrics on Lie groups \cite{MUL}. 
The Lie algebra of prederivations $\Pder(\Lg)$ forms a subalgebra of 
$\Lg\Ll (\Lg)$ containing the algebra $\Der (\Lg)$. Note that
a prederivation is just a derivation of the {\it Lie triple system}
induced by $\Lg$. A Lie triple system over $k$ is a vector space $\CT$
with a trilinear mapping $(x,y,z)\mapsto [x,y,z]$ satisfying the following
axioms \cite{HOP}:
\begin{gather*}
[x,y,z] = -[x,z,y]\\
[x,y,z]+[y,z,x]+[z,x,y]=0\\
[[x,y,z],a,b]-[[x,a,b],y,z]=[x,[y,a,b]]+[x,y,[z,a,b]]
\end{gather*}  
Any Lie algebra is at the same time a Lie triple system via $[x,y,z]:=[x,[y,z]]$.
As before in the derivation case there exist nilpotent Lie algebras 
possessing only nilpotent prederivations. In analogy to characteristically 
nilpotent Lie algebras we call a nilpotent Lie algebra {\it strongly
nilpotent} if all its prederivations are nilpotent.\\
We classify strongly nilpotent Lie algebras in dimension $7$ and 
filiform Lie algebras of dimension $n\le 11$ admitting a 
non-singular prederivation but no non-singular derivation. 
The existence of a non-singular prederivation is useful for the 
construction of affine structures \cite{BU1}
on the Lie algebra.

\section{Prederivations}

Let $\Lg$ be a Lie algebra over a field $k$ and $\Aut (\Lg)$
its automorphism group. Define a pre-automorphism of $\Lg$ to be
a bijective linear map $A : \Lg \ra \Lg$ satisfying 
\begin{equation*}
A([x,[y,z]])=[A(x),[A(y),[A(z)]]]
\end{equation*}
for all $x,y,z\in \Lg$.
The set of all pre-automorphisms forms a subgroup $\CM(\Lg)$
of $GL(\Lg)$. This group plays an important role in the study of
Lie groups which are endowed with a bi-invariant pseudo-Riemannian
metric \cite{MUL}. Over the real numbers $\CM(\Lg)$ is a closed subgroup
of $GL(\Lg)$. Its Lie algebra is denoted by $\Pder(\Lg)$ and consists
of so called pre\-derivations. One may extend this definition for any
field $k$. 

\begin{defi}
A linear map $P: \Lg \ra \Lg$ is called a prederivation of $\Lg$ if
$$P([x,[y,z]])=[P(x),[y,z]]+[x,[P(y),z]]+[x,[y,P(z)]]$$
for every $x,y,z\in \Lg$.
\end{defi}

The set of all prederivations of $\Lg$ forms a subalgebra 
$\Pder(\Lg)$ of the Lie algebra $\Lg\Ll(\Lg)$ containing the Lie algebra 
of derivations $\Der(\Lg)$: 

\begin{lem}
It holds $\Der(\Lg)\subseteq \Pder(\Lg)$.
\end{lem}

\begin{proof}
Let $D\in\Der(\Lg)$. Then by definition
$$D([x,[y,z]])=[x,D([y,z])]+[D(x),[y,z]]$$
Substituting $D([y,z])=[D(y),z]+[y,D(z)]$ we obtain
$$D([x,[y,z]])=[x,[D(y),z]]+[x,[y,D(z)]]+[D(x),[y,z]].$$
\end{proof}

Clearly we have equality for abelian Lie algebras. This is also
known to be true for semisimple Lie algebras over a field $k$ 
of characteristic zero \cite{MUL}.

\begin{prop}
Every prederivation of a finite-dimensional semisimple Lie algebra over $k$
is a derivation and hence an inner derivation:
$\Der(\Lg)=\Pder(\Lg)=\ad(\Lg)$.
\end{prop}

In the case of solvable Lie algebras equality does not hold in general.
In this paper we are mostly interested in nilpotent Lie algebras.

\begin{defi}
Let $\Lg$ be a Lie algebra and $\{ \Lg^k\}$ its lower central series
defined by $\Lg^0=\Lg,\;\Lg^k=[\Lg^{k-1},\Lg]$ for $k\ge 1$.
Recall that $\Lg$ is said to be nilpotent of degree $p$, or nilindex $p$, if
there exists an integer $p$ such that $\Lg^p=0$ and $\Lg^{p-1}\neq 0$.
A nilpotent Lie algebra of dimension $n$ and nilindex $p=n-1$ is 
called {\it filiform}.
\end{defi} 

Let $\Lg$ be a $n$-dimensional nilpotent Lie algebra. 
The descending central series
\begin{equation*}
\Lg=\Lg^0\supseteq \Lg^1 \supseteq \Lg^2 \supseteq 
\ldots \supseteq \Lg^p=(0)
\end{equation*}
defines a positive filtration $\CF$ of $\Lg$. For each $k=0,1,\ldots, p-1$ choose
a linear subspace $V_k$ so that
$$\Lg^{k}=V_{k+1}\oplus \Lg^{k+1}$$
Then as a vector space
\begin{equation}\label{filt}
\Lg=V_1\oplus V_2 \oplus \cdots \oplus V_p
\end{equation}

\begin{prop}
Let $\Lg$ nilpotent of degree $p>1$ and of dimension $n\ge 3$. 
Then $\dim \Pder (\Lg)> \dim \Der(\Lg)$.
\end{prop}

\begin{proof}
It is well known that $\dim V_1=\dim (\Lg/[\Lg,\Lg])\ge 2$ for a nilpotent
Lie algebra. Here $V_1$ denotes the subspace defined by the above filtration
$\CF$, i.e., with $\Lg=V_1\oplus [\Lg,\Lg]$. Chose a basis $(e_1,\ldots,e_r)$
of $V_1$, $r\ge 2$. Since $\Lg^2$ is strictly contained in $\Lg^1$, there 
exists a commutator of elements in $V_1$ which is not contained in $\Lg^2$. 
Hence we may assume 
$$[e_1,e_2]=e_j$$
such that $e_j\in V_2$ and $(e_1,\ldots,e_r,e_{r+1},\ldots,e_j,e_{j+1},\ldots,e_n)$
is a basis of $\Lg$ with $\Lg^2=$ $\s\{e_{j+1},\ldots,e_n\}$.
Let $Z(\Lg)$ denote the center of $\Lg$. Since $\Lg$ is nilpotent we may
choose a non-zero $z\in Z(\Lg)$. Define a linear map $P:\Lg\ra \Lg$ by
\begin{align*}
P(e_i)& =\begin{cases}
z & \text{if $i=j$},\\
0 & \text{otherwise}.
\end{cases}
\end{align*}    
Then $P$ is a prederivation of $\Lg$. Indeed, $P([a,[b,c]])=0$ for all
$a,b,c\in\Lg$ by construction, and $[P(a),[b,c]]=[a,[P(b),c]]=[a,[b,P(c)]]=0$,
since $P(\Lg)\subset Z(\Lg)$. On the other hand, assume that this $P$ 
would be a derivation. Then
$$z=P(e_j)=P([e_1,e_2])=[P(e_1),e_2]+[e_1,P(e_2)]=0$$
which is a contradiction. 
\end{proof}

\begin{lem}\label{25}
For $P\in \Pder(\Lg)$ let $\om (x,y)=P([x,y])-[x,P(y)]+[y,P(x)]$.
Then $$[x,\om(y,z)]+\om(x,[y,z])=0.$$
\end{lem}

\begin{proof}
Note that $\om \in B^2(\Lg,\Lg)$ is a $2$-coboundary for the Lie
algebra cohomology with the adjoint module. We have
\begin{align*}
[x,\om(y,z)] & = [x,P([y,z])-[x,[y,P(z)]]+[x,[z,P(y)]]\\
\om (x,[y,z]) & = P([x,[y,z]])-[x,P([y,z])]+[[y,z],P(x)]
\end{align*}
and hence, since $P$ is a prederivation
\begin{align*}
[x,\om(y,z)]+\om (x,[y,z]) & = P([x,[y,z]])-[P(x),[y,z]]-
[x,[P(y),z]]-[x,[y,P(z)]] \\
& = 0
\end{align*}
\end{proof} 

Let us now study Lie algebras admitting a non-singular prederivation.
\begin{rem}
There are many reasons why such algebras are interesting. One is the
fact that a non-singular derivation $D$ implies an {\it affine structure}
on the Lie algebra via the representation 
$\theta (x)=D^{-1}\circ \ad (x) \circ D$. The existence of affine
structures is a difficult problem with an interesting history.
For details see \cite{BU1}. If there is no such derivation,
it is still useful to have a non-singular prederivation.
The idea is to construct a bilinear product $x\cdot y:=\theta(x)y$ 
on the Lie algebra by
$$\theta (x)=P^{-1}\circ \ad (x) \circ P+\frac{1}{2}P^{-1}\circ\om(x)$$
where $\om(x)y=\om(x,y)$ as above. 
If $P\in \Der(\Lg)$ then $\om(x)=0$ and we are back to the classical
construction. In general
we have always $x\cdot y-y\cdot x=[x,y]$, but
$\theta$ might not be a representation. However, in many cases this
construction with non-singular prederivations yields affine
structures. In proposition $\ref{dim7}$ we will study $7$-dimensional
Lie algebras having only singular derivations. Many of them have
a non-singular prederivation. For example, consider
\begin{equation*}
\begin{split}
\Lg_{7,5} & = <e_1,\ldots,e_7 \mid [e_1,e_i]=e_{i+1}, i=2,3,6
\; , \; [e_1,e_4]=e_6+e_7, [e_2,e_3]=e_5,\\
& \hskip0.64 cm  [e_2,e_5]=e_6, [e_3,e_5]=e_7>
\end{split}
\end{equation*} 
Then it is easy to verify that the prederivation
$P$ defined by $P(e_1)=e_1,P(e_2)=e_2, P(e_3)=2e_3+e_4, P(e_4)=3e_4,
P(e_5)=3e_5,P(e_6)=4e_6,P(e_7)=5e_7$ induces an affine structure.
\end{rem}
 
The following generalization of Jacobson's theorem in \cite{JAC}
follows easily from Theorem $1$ in \cite{BAJ} by using the same arguments 
over $k\otimes_{\R}\C$: 

\begin{prop}
Let $\Lg$ be a Lie algebra over a field $k$ of characteristic zero admitting
a non-singular prederivation. Then $\Lg$ is nilpotent.
\end{prop}

In analogy with characteristically nilpotent Lie algebras 
one might ask whether there exist nilpotent Lie algebras possessing 
only singular
prederivations. This is indeed the case. We define a subclass of
characteristically nilpotent Lie algebras as follows:

\begin{defi}
A Lie algebra $\Lg$ over $k$ is called {\it strongly nilpotent} 
if all its prederivations are nilpotent. 
\end{defi}

\begin{rem}
Any strongly nilpotent Lie algebra is characteristically nilpotent,
but the converse is not true in general: consider the
following $7$-dimensional Lie algebra, defined by
\begin{align*}
[e_1,e_i] & = e_{i+1}, \; 2\le i\le 6 \\
[e_2,e_3] & = e_6+e_7 \\
[e_2,e_4] & = e_7 \\
\end{align*}
It is not difficult to see that all derivations are nilpotent (see also \cite{MAG}).
On the other hand, $P=\diag(1,3,3,5,5,7,7)$ is a non-singular prederivation.
One has $\dim \Pder(\Lg)=16$ and $\dim \Der(\Lg)=11$. 
\end{rem}

The following example presents
Lie algebras possessing only nilpotent prederivations.

\begin{prop}
Let $\Lg$ be the $n$-dimensional Lie algebra with basis 
$(e_1,\ldots,e_n)$, $n\ge 7$ and defining brackets
\begin{align*}
[e_1,e_i] & = e_{i+1}, \; 2\le i\le n-1 \\
[e_2,e_3] & = e_{n-1} \\
[e_2,e_4] & = e_n \\
[e_2,e_5] & = -e_n \\
[e_3,e_4] & = e_n
\end{align*}
Then $\Lg$ is strongly nilpotent.
\end{prop}

\begin{proof}
Let $P\in \Pder (\Lg)$ and write
\begin{align*}
P(e_1) & =\sum_{i=1}^n\al_ie_i\\
P(e_2) & =\sum_{i=1}^n\be_ie_i\\
P(e_3) & =\sum_{i=1}^n\ga_ie_i
\end{align*}
Note that $[e_2,e_{n-1}]=0$ since $n\ge 7$.
Using the identity
$$P([e_i,[e_j,e_k]])=[P(e_i),[e_j,e_k]]+[e_i,[P(e_j),e_k]]+[e_i,[e_j,P(e_k)]]$$
for various $(i,j,k)$ we will obtain that the associated matrix of $P$
is strictly lower-triangular.
For $(i,j,k)=(1,1,2),(1,2,4),(2,1,4),(1,2,3),(2,1,3),(3,1,3)$ and $(2,1,2)$ we obtain
\begin{align*}
P(e_4) & = (2\al_1+\be_2)e_4 +\be_3e_5+\ldots +\be_{n-4}e_{n-2}
       + \de_1e_{n-1} + \de_2e_n\\
 0     & = \be_1e_6 \\
P(e_n) & = (3\al_1+2\be_2)e_n\\
P(e_n) & = -\ga_1e_4+(\al_1+\be_2+\ga_3)e_n\\
P(e_n) & = \ga_2e_{n-1}+(3\al_1+2\be_2+\be_3-\ga_4)e_n\\
P(e_n) & = (5\al_1+2\be_2)e_n\\
P(e_{n-1}) & = (\al_1+2\be_2)e_{n-1}+(\be_3-2\be_4)e_n 
\end{align*}
It follows $\be_1=\ga_1=\ga_2=\al_1=0$ and $\ga_3=\be_2$.
Now $(i,j,k)=(1,1,k)$ for $k=3,\ldots,n-2$ yields successively
\begin{align*}
P(e_5) & = \be_2e_5+\ga_4e_6+\ldots +(2\al_2+\al_3+\ga_{n-2})e_n \\
P(e_6) & = \be_2e_6+\be_3 e_7+\ldots + (\be_{n-4}-\al_2)e_n\\
P(e_7) & = \be_2e_7+\ga_4e_8+\ldots + \ga_{n-4}e_n \\
\vdots \hskip0.3 cm & = \hskip0.3 cm \vdots \\
P(e_n) & = \be_2e_n 
\end{align*}
Comparing with $P(e_n)=2\be_2e_n$ we obtain $\be_2=0$ and $P$ is nilpotent.
\end{proof}

\begin{example}
Let $n=7$ and $\Lg$ as above. Then we see that
the algebra $\Pder(\Lg)$ is given by the set of the following matrices:
\begin{equation*}
\begin{pmatrix} 
0 & 0 & 0 & 0 & 0 & 0 & 0\\
\al_2 & 0 & 0 & 0 & 0 & 0 & 0\\
\al_3 & \be_3 & 0 & 0 & 0 & 0 & 0\\
\al_4 & \frac{\al_2}{2} & \be_3 & 0 & 0 & 0 & 0\\
\al_5 & \be_5 & \ga_5 & \be_3 & 0 & 0 & 0 \\
\al_6 & \be_6 & \ga_6 & \frac{3\al_2}{2} & \be_3 & 0 & 0\\
\al_7 & \be_7 & \ga_7 & \be_5-\al_3-\al_4 & \ga_5+2\al_2+\al_3 &
\be_3 -\al_2 & 0\\
\end{pmatrix}
\end{equation*} 
One has $\dim \Pder(\Lg)=13$ and $\dim \Der(\Lg)=10$.
\end{example}

We can also construct a series of characteristically
nilpotent Lie algebras which are not strongly nilpotent.
Let $\Lg_n(\al)$ for $\al\ne 0$ be the Lie algebras with 
basis $(e_1,\ldots,e_n)$, $n\ge 7$ and defining brackets
\begin{align*}
[e_1,e_i] & = e_{i+1}, \; 2\le i\le n-1 \\
[e_2,e_3] & = e_5+\al e_n \\
[e_2,e_j] & = e_{j+2},\; 4\le j\le n-2 \\
\end{align*}
Over the complex numbers, $\Lg_n(\al)\cong \Lg_n(1)=\Lg_n$. An isomorphism
is given by $\phi (e_i)=\la^i e_i$ with a $\la$ satisfying
$\la^{n-5}=\al^{-1}$. 

\begin{prop}
For every $n\ge 7$ the Lie algebra $\Lg_n$ is 
characteristically nilpotent, but not strongly nilpotent.
\end{prop}

\begin{proof}
We construct a non-singular prederivation of $\Lg_n$ by
\begin{align*}
P(e_i)&=\begin{cases}
3e_i+(5-n)e_{n-2} & \text{if $i=3$},\\
5e_i+(5-n)e_{n} & \text{if $i=5$},\\
ie_i & \text{otherwise}.
\end{cases}
\end{align*} 
Then $\det P=n!\ne 0.$ We have to ckeck the identity
$$P([e_i,[e_j,e_k]])=[P(e_i),[e_j,e_k]]+[e_i,[P(e_j),e_k]]+[e_i,[e_j,P(e_k)]]$$
for $1\le i,j,k \le n$, where we may assume that $j<k$ and $i\le k$. 
The identity clearly holds in all cases where $P(e_3)$ or $P(e_5)$ is
not involved since then $P(e_i)=ie_i$ and hence 
$P([e_i,e_j])=[P(e_i),e_j]+[e_i,P(e_j)]$ for $i,j,i+j\ne 3,5$.
The term $[e_i,[e_j,e_k]]$ belongs to $\s \{e_4,\ldots e_n\}$. It equals 
$e_5$ for $(i,j,k)=(1,1,3)$ or $(i,j,k)=(2,1,2)$ in which case the above 
identity reads as $P(e_5)=2e_5+3e_5+(5-n)e_n$, respectively 
$P(e_5+e_n)=(2+1+2)(e_5+e_n)$.
It remains to ckeck the cases where $i,j$ or $k$ equals $3$ or $5$.
But this is easily done. We have $[e_i,P(e_5)]=5[e_i,e_5]$ since 
$e_n$ belongs to the center of $\Lg_n$. It holds $[e_i,P(e_3)]=0$ for 
$i\ge 3$ and $[e_1,P(e_3)]=3e_4+(5-n)e_{n-1}$, $[e_2,P(e_3)]=3e_5+3e_n$. \\
Now assume that $D\in \Der(\Lg_n)$ is a derivation. Write
\begin{align*}
D(e_1) & =\sum_{i=1}^n\ze_ie_i\\
D(e_2) & =\sum_{i=1}^n\mu_ie_i
\end{align*}
Using $D(e_i)=[D(e_1),e_{i-1}]+[e_1,D(e_{i-1})]$ we  
compute $D(e_3),D(e_4),\ldots D(e_n)$ successively.
We have 
\begin{align*}
D(e_3) & = (\ze_1+\mu_2)e_3+\mu_3e_4+(\mu_4-\ze_3)e_5+\ldots+
(\mu_{n-1}-\ze_{n-2}-\ze_3)e_n\\
D(e_4) & = (2\ze_1+\mu_2)e_4+(\ze_2+\mu_3)e_5 + \ldots + 
(\ze_2+\mu_{n-2}-\ze_{n-3})e_n\\
D(e_5) & =(3\ze_1+\mu_2)e_5+(\mu_3+2\ze_2)e_6+\ldots+(\mu_{n-3}-\ze_{n-4})e_n\\ 
\vdots \hskip0.3 cm & = \hskip0.3 cm \vdots\\
D(e_n) & =((n-2)\ze_1+\mu_2)e_n\\
\end{align*} 
Then $D([e_2,e_3])=[D(e_2),e_3]+[e_2,D(e_3)]$ is, for $n\ge 7$,
equivalent to
\begin{equation*}
-\mu_1e_4+(2\ze_1-\mu_2)e_5+2\ze_2e_6+((n-3)\ze_1-\mu_2)e_n = 0,
\end{equation*} 
which implies $\mu_1=0,\, \mu_2=2\ze_1,\, (n-5)\ze_1=0$ and hence 
$\mu_1=\mu_2=\ze_1=0$. It follows that the matrix for $D$
is strictly lower triangular, hence $D$ is nilpotent.
Note that we have used $n\ge 7$ in the computation. 
For $n\le 6$ there exists always a non-singular derivation of $\Lg_n$.
\end{proof}

Let $\Lg$ a $p$-step nilpotent Lie algebra. If $p<3$ then $\Lg$ always admits 
a non-singular derivation. However there are examples of $3$-step nilpotent
Lie algebras possessing only nilpotent derivations \cite{DIL}.
The result for prederivations is as follows.

\begin{prop}
For $p<5$ any $p$-step nilpotent Lie algebra admits a non-singular
prederivation.
\end{prop}

\begin{proof}
For a $2$-step nilpotent or abelian Lie algebra, every linear
map $P:\Lg\ra \Lg$ is a prederivation. Then the claim is obvious.
Hence we may assume that $p=3$ or $p=4$, and $\Lg^4=0$.  
Then \eqref{filt} says
$$\Lg=V_1\oplus V_2\oplus V_3\oplus V_4$$ 
Let $(e_1,\ldots,e_n)$ be a basis of $\Lg$ adapted to this decomposition,
i.e., which is the union of the bases for $V_i$. Define a linear map 
$P:\Lg \ra \Lg$ by
\begin{align*}
P(e_i)&=\begin{cases}
e_i  & \text{if $e_i \in V_1$ or $V_2$},\\
3e_i & \text{if $e_i \in V_3$ or $V_4$.}\
\end{cases}
\end{align*} 
In particular $P(e_i)=3e_i$ for all $e_i\in \Lg^2$. Writing $P(e_i)=\ze_i e_i$
for $i=1,\ldots,n$, $P$ is a prederivation if and only if
$$P([e_i,[e_j,e_k]])=(\ze_i+\ze_j+\ze_k)[e_i,[e_j,e_k]]$$
for all $i,j,k$. The term $[e_i,[e_j,e_k]]$ is zero for $e_j\in \Lg^2$ or
$e_k\in \Lg^2$ because of $\Lg^4=0$. If $e_i\in  \Lg^2$, then we apply the
Jacobi identity to see that
$$[e_i,[e_j,e_k]]=[e_j,[e_i,e_k]]+[e_k,[e_i,e_j]]$$
is contained in $\Lg^4$ which is zero. Hence the above identity for $P$
is trivially satisfied if not $e_i,e_j,e_k \in V_1$ or $V_2$. But then we have
$[e_i,[e_j,e_k]]\in \Lg^2$ and the identity is equivalent to
$P([e_i,[e_j,e_k]])=3 [e_i,[e_j,e_k]]$. Hence $P$ is a non-singular prederivation.
\end{proof}

How big is the class of strongly nilpotent Lie algebras ?
In dimension $7$ there are already infinitely many complex non-isomorphic
strongly nilpotent Lie algebras:

\begin{prop}\label{dim7}
Any strongly nilpotent complex Lie algebra of dimension $7$ is isomorphic
to one of the following algebras:
\begin{equation*}
\begin{split}
\Lg_{7,1} & = <e_1,\ldots,e_7 \mid [e_1,e_i]=e_{i+1}, 2\le i\le 6,
\; [e_2,e_3]=e_6, [e_2,e_4]=e_7,\\
& \hskip0.64 cm  [e_2,e_5]=e_7,[e_3,e_4]=-e_7>\\
\Lg_{7,4}^{\la} & =<e_1,\ldots,e_7 \mid [e_1,e_2]=e_3, [e_1,e_3]=e_4,
[e_1,e_4]=e_6+\la e_7,[e_1,e_5]=e_7,\\
& \hskip0.64 cm  [e_1,e_6]=e_7,[e_2,e_3]=e_5,[e_2,e_4]=e_7,[e_2,e_5]=e_6,
[e_3,e_5]=e_7>\\
\Lg_{7,7} & = <e_1,\ldots,e_7 \mid [e_1,e_i]=e_{i+1}, i=2,3,6,\; [e_1,e_4]=e_7,
[e_1,e_5]=e_7,\\
& \hskip0.64 cm [e_2,e_3]=e_5,[e_2,e_4]=e_7,[e_2,e_5]=e_6,[e_3,e_5]=e_7>
\end{split}
\end{equation*} 
These algebras are pairwise non-isomorphic except for 
$\Lg_{7,4}^{\la}\simeq \Lg_{7,4}^{-\la}$.
\end{prop}
\begin{proof}
Since a characteristically nilpotent complex Lie algebra of
dimension $7$ is indecomposable, we may take the classification
of such algebras from \cite{MAG}: the list contains the algebras
$\Lg_{7,1},\ldots,\Lg_{7,8}$, where $\Lg_{7,4}$ depends on a 
parameter $\la\in\C$. Now one has to compute the algebra of
prederivations in each case. We did this, but it does not seem 
useful to write down all the computations here.  
The following table shows the results:
\begin{center}
\begin{tabular}{|c|c|c|c|}
\hline
Algebra & $\dim \Der (\Lg)$ & $\dim \Pder (\Lg)$  & $P^{-1}$ exists \\
\hline\hline
$\Lg_{7,1}$ & $10$ & $13$ & $\mi$   \\ \hline
$\Lg_{7,2}$ & $10$ & $13$ & $\checkmark$   \\ \hline
$\Lg_{7,3}$ & $11$ & $16$ & $\checkmark$          \\ \hline
$\Lg_{7,4}$ & $10$ & $12$ & $\mi$   \\ \hline
$\Lg_{7,5}$ & $10$ & $13$ & $\checkmark$   \\ \hline
$\Lg_{7,6}$ & $10$ & $13$ & $\checkmark$  \\ \hline
$\Lg_{7,7}$ & $10$ & $12$ & $\mi$         \\ \hline
$\Lg_{7,8}$ & $10$ & $14$ & $\checkmark$   \\ \hline
\end{tabular} 
\end{center}      
The algebras without a non-singular prederivation are in fact
strongly nilpotent.
\end{proof}

\section{Prederivations of filiform Lie algebras}

If $\Lg$ is a filiform Lie algebra of dimension $n$,
then there exists an adapted basis $(e_1,\ldots,e_n)$ for $\Lg$, see
\cite{BU1}. 
The brackets of such a filiform Lie algebra
with respect to the basis $(e_1,\ldots,e_n)$ are then given by

\begin{align}\label{lie}
[e_1,e_i] & =e_{i+1}, \quad i=2,\dots ,n-1 \\
[e_i,e_j] & =\sum_{r=1}^n\biggl(\;\sum_{\ell=0}^{[(j-i-1)/2]} (-1)^\ell
\binom{j-i-\ell-1}{\ell}\al_{i+\ell,\, r-j+i+2\ell+1}\biggr)e_r,
 \quad 2 \le i<j \le n.
\end{align}
with constants $\al_{k,s}$ which are zero for all pairs $(k,s)$ not in 
the index set $\CI_n$. Here 
$\CI_n$ is given by
\begin{align*}
\CI_n^0 &=\{(k,s)\in \N \times \N \mid 2 \le k \le [n/2],\,
2k+1 \le s \le n \},\\
\CI_n& =\begin{cases}
\CI_n^0 & \text{if $n$ is odd},\\
\CI_n^0 \cup \{(\frac{n}{2},n)\} & \text{if $n$ is even}.
\end{cases}
\end{align*}     

For filiform Lie algebras we study the conditions for
the existence of a non-singular prederivation. As it turns
out there are only a few algebras possessing a non-singular 
prederivation but no non-singular derivation:

\begin{prop}
Up to isomorphism there are the following  
filiform Lie algebras $\Lg$ of dimension $n\le 11$ over $\C$ 
possessing a non-singular prederivation but no non-singular
derivation:

\begin{center}
\begin{tabular}{|c|c|}
\hline
$\dim \Lg$ & Algebra \\
\hline\hline
$7$ & $\mu_7^4,\mu_7^6$  \\ \hline
$8$ & $\mu_8^{11}(0),\mu_8^{13}$   \\ \hline
$9$ & $\mu_9^{28}(0),\mu_9^{30},\mu_9^{32},\mu_9^{34}$ \\ \hline
$10$ & $\mu_{10}^8(0,0,0),\mu_{10}^{15}(0,0),\mu_{10}^{48}$   \\ \hline
$11$ & $\mu_{11}^4(0,0),\mu_{11}^{16},\mu_{11}^{71},\mu_{11}^{62}
(0,\be,0,0),\be \ne 0$, \\ \hline
     & $\mu_{11}^{89}, \mu_{11}^{81}, \mu_{11}^{96a}, \mu_{11}^{101a}$   \\ \hline
\end{tabular} 

\end{center}  
Here we use the notation from the classification list of \cite{GJK}.
\end{prop}

\begin{rem}
We have found $2$ new filiform Lie algebras $\mu_{11}^{96a}, \mu_{11}^{101a}$
given by 
\begin{align*}
\mu_{11}^{96a} & = \mu_0+\Psi_{1,8}+\Psi_{1,10}\\
\mu_{11}^{101a} & = \mu_0+\Psi_{1,9}+\Psi_{1,10}
\end{align*}
It seems that they are not isomorphic to one of the algebras in the list
of \cite{GJK}.
\end{rem}

\begin{proof}
Since the filiform Lie algebras in question are classified, it would 
be possible to prove the result by calculating the derivations and prederivations
separately for each algebra of the classification list.
We proceed differently, however. Since we may write any
filiform Lie algebra with respect to an adapted basis as in \eqref{lie} and $(3)$, 
we can determine the algebras possessing a 
non-singular prederivation, but no non-singular derivation,
with respect to the structure constants $\al_{i,j}$.
Then we obtain a small list of algebras not necessarily beeing
non-isomorphic. The result follows then by determining 
the isomorphisms between the remaining algebras. In all but two cases 
(see above) we could easily find an isomorphism to an algebra of the list
in \cite{GJK}.\\
We present the computations for $\dim \Lg=11$. The Lie brackets
relative to an adapted basis $(e_1,\ldots,e_{11})$ are given 
by \eqref{lie} and $(3)$: \\
\begin{align*}
\begin{split}
[e_1,e_i] & =e_{i+1}, \; 2\le i\le 9 \\
[e_2,e_3] & = \al_{2,5}e_5+\al_{2,6}e_6+\al_{2,7}e_7 +\al_{2,8}e_8+
\al_{2,9}e_9+\al_{2,10}e_{10}+\al_{2,11}e_{11}\\
[e_2,e_4] & = \al_{2,5}e_6+\al_{2,6}e_7+\al_{2,7}e_8+\al_{2,8}e_9 +
\al_{2,9}e_{10}+\al_{2,10}e_{11}\\
[e_2,e_5] & = (\al_{2,5}-\al_{3,7})e_7+ (\al_{2,6}-\al_{3,8})e_8+
(\al_{2,7}-\al_{3,9})e_9\\
& \quad +(\al_{2,8}-\al_{3,10})e_{10}+ (\al_{2,9}-\al_{3,11})e_{11}\\
[e_2,e_6] & = (\al_{2,5}-2\al_{3,7})e_8+ (\al_{2,6}-2\al_{3,8})e_9 +
(\al_{2,7}-2\al_{3,9})e_{10}+(\al_{2,8}-2\al_{3,10})e_{11}\\
[e_2,e_7] & = (\al_{2,5}-3\al_{3,7}+\al_{4,9})e_9 +(\al_{2,6}-3\al_{3,8}+
\al_{4,10})e_{10}\\
& \quad +(\al_{2,7}-3\al_{3,9}+\al_{4,11})e_{11}\\
[e_2,e_8] & = (\al_{2,5} - 4\al_{3,7} + 3\al_{4,9})e_{10}+
(\al_{2,6} - 4\al_{3,8} + 3\al_{4,10})e_{11}\\
[e_2,e_9] & = (\al_{2,5}-5\al_{3,7}+6\al_{4,9}-\al_{5,11})e_{11}\\
%\end{split}
%\end{align*}  
%\begin{align*}
%\begin{split}
[e_3,e_4] & = \al_{3,7}e_7+\al_{3,8}e_8+\al_{3,9}e_9+\al_{3,10}e_{10}+
\al_{3,11}e_{11}\\
[e_3,e_5] & = \al_{3,7}e_8+\al_{3,8}e_9+\al_{3,9}e_{10}+\al_{3,10}e_{11} \\
[e_3,e_6] & = (\al_{3,7}-\al_{4,9})e_9 +(\al_{3,8}-\al_{4,10})e_{10}+
(\al_{3,9}-\al_{4,11})e_{11}\\
[e_3,e_7] & = (\al_{3,7}-2\al_{4,9})e_{10}+(\al_{3,8}-2\al_{4,10})e_{11}\\
[e_3,e_8] & = (\al_{3,7}-3\al_{4,9}+\al_{5,11})e_{11}
\end{split}
\end{align*}  
\begin{align*}
\begin{split}
[e_4,e_5] & = \al_{4,9}e_9 + \al_{4,10}e_{10}+\al_{4,11}e_{11}\\
[e_4,e_6] & = \al_{4,9}e_{10}+\al_{4,10}e_{11}\\
[e_4,e_7] & = (\al_{4,9}-\al_{5,11})e_{11}\\
[e_5,e_6] & = \al_{5,11}e_{11}\\
\end{split}
\end{align*}    
The Jacobi identity is satisfied if and only if:
\begin{equation*}
\begin{split}
0 & =\al_{4,9}(2\al_{2,5}+\al_{3,7}) - 3\al_{3,7}^2\\
0 & =\al_{4,10}(2\al_{2,5}+\al_{3,7})+3\al_{4,9}(\al_{2,6}+\al_{3,8})
-7\al_{3,7}\al_{3,8}\\
0 &= \al_{5,11}(2\al_{2,5}-\al_{3,7}-\al_{4,9})+\al_{4,9}(6\al_{4,9}-4
\al_{3,7})\\
0 & =\al_{4,12}(2\al_{2,7}+\al_{3,9})+\al_{4,11}(2\al_{2,5}+\al_{3,7})+
3\al_{4,10}(\al_{2,6}+\al_{3,8})-4\al_{3,8}^2\\
  & \; + 2\al_{4,9}(2\al_{2,7}+3\al_{3,9})-8\al_{3,7}\al_{3,9}
\end{split}
\end{equation*}  
We will denote the fact that there exists a non-singular
prederivation $P\in \Pder (\Lg)$ simply by "$P^{-1}$ exists".
We have to distinguish several cases.
Consider first the case $$2\al_{2,5}+\al_{3,7}=0$$
If $\al_{4,9}\ne 0$ then $P^{-1}$ exists if and only if
\begin{align*}
\al_{2,6} & = 0\\
\al_{2,8} & = 0\\
\al_{3,11} & = (\al_{3,10}\al_{4,10}+4\al_{4,9}\al_{2,9}+\al_{3,9}^2)/\al_{4,9}\\
\al_{4,11} & = (\al_{4,10}^2+6\al_{4,9}\al_{3,9})/\al_{4,9}
\end{align*}
A non-singular derivation $D^{-1}$ exists if and only if the above
conditions are satisfied. Hence in this case "$P^{-1}$ exists if and only
if $D^{-1}$ exists". \\
In the following let $\al_{4,9}=0$. 
For $\al_{2,6},\al_{3,8}\ne 0$ we have: $P^{-1}$ 
exists if and only if $D^{-1}$ exists.\\
For $\al_{2,6}\ne 0$, $\al_{3,8}=0$ $P^{-1}$ exists if and only if
\begin{align*}
\al_{2,i} & = 0, \; i=7,9,11\\
\al_{2,10} & = 4\al_{2,8}^2/(3\al_{2,6})\\
\al_{3,i} & =0, \; i=7,8,9,10,11\\
\al_{4,i} & =0, \; i=9,10,11\\
\al_{5,11}& =0
\end{align*}
However $D^{-1}$ exists if and only if the above conditions hold
{\it and} $\al_{2,11}=0$.
Hence we have the following algebra, given by $16$ structure constants
$(\al_{2,5},\al_{2,6},\ldots,\al_{5,11})$ as follows:
$$(0,a_2,0,a_4,0,\frac{4a_4^2}{3a_2},a_7,0,0,0,0,0,0,0,0,0)$$
where $a_2,a_7\ne 0$. This algebra is isomorphic over $\C$ to
$\mu_{11}^{71}$ given by
$$(0,1,0,0,0,0,1,0,0,0,0,0,0,0,0,0)$$
In the following let $\al_{2,6}= 0$. Consider the case $\al_{2,7},\al_{3,9}\ne 0$.
If $\al_{4,11}\ne 2\al_{3,9}$ then $P^{-1}$ exists if and only if
$D^{-1}$ exists. If $\al_{4,11}= 2\al_{3,9}$ then $P^{-1}$ exists if 
and only if
\begin{align*}
\al_{2,8} & = 0\\
\al_{3,10} & = 0\\
\al_{3,11} & =\al_{4,11}\al_{2,9}/\al_{3,9}\\
\al_{4,10} & =0\\
\al_{5,11}& =0
\end{align*}
On the other hand $D^{-1}$ exists if and only if in addition
$\al_{2,11}=\al_{2,9}^2/\al_{3,9}$. We obtain the algebra
$$(0,0,a_3,0,a_5,a_6,a_7,0,0,a_{10},0,2a_5,0,0,2a_{10},0)$$
where $a_3\ne 0, a_7\ne a_5^2/a_{10}$. It is isomorphic over $\C$ to
$\mu_{11}^{62}(0,\be,0,0), \be=a_3/(2a_{10})\ne 0$ given by
$$(0,0,\be,0,0,0,1,0,0,\frac{1}{2},0,0,0,0,1,0)$$ 
Consider the case $\al_{2,7}\ne 0,\al_{3,9}=0$. If $\al_{4,11}\ne 0$ 
then $P^{-1}$ exists if and only if $D^{-1}$ exists.
If  $\al_{4,11}= 0$ then $P^{-1}$ exists if and only if
\begin{align*}
\al_{2,i} & = 0, \; i=8,9\\
\al_{3,i} & = 0, \; i=7,8,10,11\\
\al_{4,10}& =0\\
\al_{5,11}& =0
\end{align*}
whereas $D^{-1}$ exists if and only in addition $\al_{2,11}=0$.
We obtain the algebra
$$(0,0,a_3,0,0,a_6,a_7,0,0,0,0,0,0,0,0,0)$$
where $a_3,a_7\ne 0$. It is isomorphic over $\C$ to
$\mu_{11}^{81}$ given by
$$(0,0,1,0,0,0,1,0,0,0,0,0,0,0,0,0)$$
In the following let $\al_{2,7}=0$. If $\al_{3,8}\ne 0$
then  $P^{-1}$ exists if and only if $D^{-1}$ exists.
This is also true for $\al_{3,8}=0, \al_{3,9}\ne 0$. 
Assume in the following $\al_{3,8}=\al_{3,9}=0$.
If $\al_{2,8}\ne 0$ then $P^{-1}$ exists if 
and only if $D^{-1}$ exists except for the case where
$\al_{3,10}=0,\al_{2,11}\ne 0$. We obtain the algebra
$$(0,0,0,a_4,0,0,a_7,0,0,0,0,0,0,0,0,0)$$
where $a_4,a_7\ne 0$. It is isomorphic over $\C$ to $\mu_{11}^{89}$
given by
$$(0,0,0,1,0,0,1,0,0,0,0,0,0,0,0,0)$$
Assume in the following $\al_{2,8}=0$. If $\al_{3,11}\ne 0$
then $P^{-1}$ exists if and only if $D^{-1}$ exists.
Assume $\al_{3,11}=0$. If $\al_{2,9}\ne 0$ then $P^{-1}$ exists 
if and only if
\begin{align*}
\al_{2,10} & = 0\\
\al_{3,10} & = 0\\
\al_{4,i}& =0, \, i=10,11\\
\al_{5,11}& =0
\end{align*}
However $D^{-1}$ exists if and only if in addition $\al_{2,11}=0$.
We obtain the algebra
$$(0,0,0,0,a_5,0,a_7,0,0,0,0,0,0,0,0,0)$$
where $a_5,a_7\ne 0$. It is isomorphic over $\C$ to
$$(0,0,0,0,1,0,1,0,0,0,0,0,0,0,0,0)$$
For $\al_{2,9}=0$ we have always that $P^{-1}$ exists if and only if 
$D^{-1}$ exists except for the algebra
$$(0,0,0,0,0,a_6,a_7,0,0,0,0,0,0,0,0,0)$$
where  $a_6,a_7\ne 0$. It is isomorphic over $\C$ to
$$(0,0,0,0,0,1,1,0,0,0,0,0,0,0,0,0)$$
To finish the proof we have to consider the case
$$2\al_{2,5}+\al_{3,7}\ne 0$$
If $\al_{2,5}=0,\al_{3,7}\ne 0$ then $P^{-1}$ exists if and only if 
$D^{-1}$ exists. If $\al_{2,5}\ne 0,\al_{3,7}=0$ then
we obtain an algebra which is isomorphic to $\mu_{11}^{16}$, given by
$$(1,0,0,0,0,0,1,0,0,0,0,0,0,0,0,0)$$
Assume in the following $\al_{2,5},\al_{3,7}\ne 0$. Then 
$P^{-1}$ exists if and only if $D^{-1}$ exists except for the case
where $7\al_{3,7}=4\al_{2,5}, 5\al_{3,7}^2-6\al_{2,5}\al_{3,7}+4\al_{2,5}
\ne 0$. In that case we obtain an algebra being isomorphic to
$\mu_{11}^4(0,0)$ given by
$$(1,0,0,0,0,0,1,\frac{4}{7},0,0,0,0,\frac{8}{21},0,0,0)$$
That concludes our proof.
\end{proof}

We have also studied the existence question of non-singular
prederivations for certain filiform algebras of dimension $n\ge 12$.
We will consider the following conditions, which are isomorphism invariants of $\Lg$:
 \begin{itemize}
\item[(a)] $\Lg$ contains no one-codimensional subspace $U \supseteq \Lg^1$
such that $[U,\Lg^1]\subseteq \Lg^4$.
\item[(b)] $\Lg^{\frac{n-4}{2}}$
is abelian, if $n$ is even.
\item[(c)] $[\Lg^1,\Lg^1]\subseteq \Lg^6$.
\end{itemize}

We will focus on algebras satifying property $(a)$:

\begin{defi}
Let $\CA_n^1$ denote the set of $n$-dimensional filiform laws whose algebras
satisfy the properties $(a),(b),(c)$. Denote by $\CA_n^2$ the set of
$n$-dimensional filiform laws whose algebras satisfy $(a),(b)$, but {\it not}
$(c)$. 
\end{defi}

The above properties of $\Lg$ can be expressed in terms of the
corresponding structure constants $\al_{k,s}$. It is easy to verify the
following (use \eqref{lie}, $(3)$):

\begin{itemize}
\item[] $\al_{2,5}\ne 0$, if and only $\Lg$ satisfies property $(a)$.
\item[] $\al_{\frac{n}{2},n}=0$, if and only if $\Lg$ satisfies
property $(b)$.
\item[] $\al_{3,7}=0$, if and only if $\Lg$ satisfies property $(c)$.
\end{itemize}

If $\Lg$ satisfies property $(a)$ we may change the adpated basis so that
it stays adapted and
$$\al_{2,5}=1.$$
In fact, we may take $f\in GL(\Lg)$ defined by $f(e_1)=ae_1, f(e_2)=be_2$
and $f(e_i)=[f(e_1),f(e_{i-1})]$ for $3\le i\le n$ with suitable
nonzero constants $a$ and $b$.\\                                               
The following results follow by straightforward computation:

\begin{prop}
Let $\Lg$ be a filiform Lie algebra with law in $\CA_n^1$. Then 
$\Lg$ admits a non-singular prederivation if and only if
\begin{align*}
\al_{3,i} & = 0, \quad i=8,\ldots,n\\
\al_{2,j} & = \frac{1}{2^{j-5}(j-3)}\binom{2j-8}{j-4}\al_{2,6}^{j-5},\quad 
j=7,\ldots,n-1
\end{align*}
On the other hand, $\Lg$ admits a non-singular derivation if and only if
in addition 
$$ \al_{2,n}=\frac{1}{2^{n-5}(n-3)}\binom{2n-8}{n-4}\al_{2,6}^{n-5}$$
\end{prop}

Note that the formula contains the Catalan numbers 
$$C_m=\frac{1}{m+1}\binom{2m}{m}$$

\begin{prop}
Let $\Lg$ be a filiform Lie algebra with law in $\CA_n^2$. Then
$\Lg$ admits a non-singular prederivation if and only if it admits
a non-singular derivation.
\end{prop}
\vskip0.5 cm
\small
{\bf Acknowledgement:} I am grateful to the referee for helpful
remarks.

\end{document}